\documentclass[twoside,11pt,a4paper]{amsart}
\usepackage{amssymb,supertabular,a4}
\usepackage{amsmath,amsthm,amscd,amssymb,graphicx}
\input xy
\xyoption{all}
\def\Hom{\hbox{\rm Hom\,}}
\def\Ext{\hbox{\rm Ext\,}}

\def\Hom{\hbox{\rm Hom\,}}
\def\dimHom{\hbox{\rm dimHom\,}}
\def\dimExt{\hbox{\rm dimExt\,}}

\def\mod{\hbox{\rm mod\,}}

\def\ndp{\hbox{\rm ndp\,}}

\def\ql{\hbox{\rm ql\,}}
\def\udim{\hbox{\rm \underline{dim}\,}}
\def\ra{{\rightarrow}}

\title[The number of GR measures processing no direct predecessors]
{The number of the Gabriel-Roiter measures admitting no direct
predecessors over a wild quiver}
\author{Bo Chen}
\date{}
\address {Universit\"at zu K\"oln\\
          Mathematisches Institut\\
                   Weyertal 86-90\\
           D-50931 K\"oln\\ Germany}

\email {mcebbchen@googlemail.com}

\thanks{The author is supported by DFG-Schwerpunktprogramm
`Representation theory'.}

\begin{document}

\maketitle

\newtheorem{theo}{Theorem}[section]
\newtheorem{defi}[theo]{Definition}
\newtheorem{lemm}[theo]{Lemma}
\newtheorem{coro}[theo]{Corollary}
\newtheorem{prop}[theo]{Proposition}
\newtheorem*{main}{Main Property}
\newtheorem{Example}{Example}
\newtheorem*{Question}{Question}

\thispagestyle{empty}

\begin{abstract} Let $Q$ be the wild quiver with three vertices,
labeled by $1$,$2$ and $3$, and one arrow from $1$ to $2$ and two
arrows from $2$ to $3$. The Gabriel-Roiter submodules of the
indecomposable preprojective modules and some quasi-simple modules
are described using the embedding of the Kronecker modules. It will
be shown that there are infinitely many Gabriel-Roiter measures
admitting no direct predecessors.
\end{abstract}

\bigskip



\section{Introduction}

Let $\Lambda$ be an artin algebra and $\mod\Lambda$ the category of
finitely generated right $\Lambda$-modules. For each
$M\in\mod\Lambda$, we denote by $|M|$ the length of $M$. The symbol
$\subset$ is used to denote proper inclusion. The Gabriel-Roiter (GR
for short) measure $\mu(M)$ for a $\Lambda$-module $M$ was defined
in \cite{R2} by induction as follows:
$$\mu(M)=\left\{\begin{array}{ll}
0, & \textrm{if $M=0$};\\
\max_{N\subset M}\{\mu(N)\}, & \textrm{if $M$ is decomposable};\\
\max_{N\subset M}\{\mu(N)\}+\frac{1}{2^{|M|}}, & \textrm{if $M$ is indecomposable}.\\
\end{array}\right.$$
(In later discussion, we will use the original definition for our
convenience, see \cite{R1} or section \ref{GRdef} below.) The
so-called Gabriel-Roiter submodules of an indecomposable module are
defined to be the indecomposable proper submodules with maximal GR
measure.

Using Gabriel-Roiter  measure, Ringel obtained a partition of the
module category for any artin algebra of infinite representation
type \cite{R1,R2}: there are infinitely many GR measures $\mu_i$ and
$\mu^i$ with $i$ natural numbers, such that
$$\mu_1<\mu_2<\mu_3<\ldots\quad \ldots <\mu^3<\mu^2<\mu^1$$ and such that any
other GR measure $\mu$ satisfies $\mu_i<\mu<\mu^j$ for all $i,j$.
The GR measures $\mu_i$ (resp. $\mu^i$) are called take-off (resp.
landing) measures. Any other GR measure is called a central measure.
An indecomposable module  is  called a  take-off (resp. central,
landing) module if its GR measure is a take-off (resp. central,
landing) measure.

Ringel showed in \cite{R1} that all modules lying in the landing
part are preinjective modules in the sense of Auslander and Smal\o\
\cite{AS}. In \cite{Ch4}, it was shown that for tame quivers, all
indecomposable preprojective modules are take-off modules.

Given an artin algebra $\Lambda$. Let $\mu,\mu'$ be two GR measures
for $\Lambda$. We call $\mu'$  a {\bf direct successor} of $\mu$ if,
first, $\mu<\mu'$ and second, there does not exist a GR measure
$\mu''$ with $\mu<\mu''<\mu'$. The so-called {\bf Successor Lemma}
in \cite{R2} states that any Gabriel-Roiter measure $\mu$ different
from $\mu^1$, the maximal one, has a direct successor. However,
there is no `Predecessor Lemma'. For example, the minimal central
measure (if exists) has no direct predecessor. We denote by
$\ndp(\Lambda)$ the number of the GR measures admitting no direct
predecessors. It is clear that any GR measure over a
representation-finite artin algebra has a direct predecessor. It was
shown in \cite{Ch6} that for each tame hereditary algebra over an
algebraically closed field, $1\leq \ndp(\Lambda)<\infty$. Thus we
may ask the following question:

\begin{Question} Does the number of the GR measures having no direct
predecessors relate to the representation type of artin algebras?
More precisely,  does a representation-infinite (hereditary) algebra
(over an algebraically closed field) is of wild type imply that
there are  infinitely many GR measures having no direct predecessors
and vice versa?
\end{Question}

However, for wild quivers, or more general wild algebras,  it is
difficult to calculate the Gabriel-Roiter measures of the
indecomposable modules or to show that whether a GR measure has a
direct predecessor or not. In this paper, we will mainly study the
wild quiver $Q$
$$\xymatrix{1\ar@{->}[r] & 2\ar@/^/[r]\ar@/_/[r]& 3\\}$$ and the
category of finite dimensional representations (simply called
modules) over an algebraically closed field $k$.  After some
preliminaries, we will calculate the GR submodules of each
indecomposable preprojective module (Proposition \ref{P}) and
describe the take-off part (Proposition \ref{takeoff}). Moreover, we
will show the existence of the minimal central measure, thus the
existence of the Gabriel-Roiter measure having no direct predecessor
(Proposition \ref{minimal}). In particular, we will see that,
different from tame quiver cases, the take-off part contains only
finitely many preprojective modules. An indecomposable module with
dimension vectors $(0,a,b)$ is called a Kronecker module. The
Gabriel-Roiter submodules of the quasi-simple modules of the form
$\tau^{-i}M,i\geq 0$ are described, where $M$ is a Kronecker module
and $\tau$ is the Auslander-Reiten translation (Proposition
\ref{H_a}, \ref{H(a)}, \ref{H^a}). Finally, we will show that there
are infinitely many GR measures admitting no direct predecessors
(Theorem \ref{ndp}).

\section{Preliminaries and examples}
\subsection{The Gabriel-Roiter measure}\label{GRdef}
We fist  recall the original definition of Gabriel-Roiter measure
\cite{R1,R2}. Let $\mathbb{N}_1$=$\{1,2,\ldots\}$ be the set of
natural numbers and $\mathcal{P}(\mathbb{N}_1)$ be the set of all
subsets of $\mathbb{N}_1$.  A total order on
$\mathcal{P}(\mathbb{N}_1)$ can be defined as follows: if $I$,$J$
are two different subsets of $\mathbb{N}_1$, write $I<J$ if the
smallest element in $(I\backslash J)\cup (J\backslash I)$ belongs to
J. Also we write $I\ll J$ provided $I\subset J$ and for all elements
$a\in I$, $b\in J\backslash I$, we have $a<b$. We say that $J$ {\bf
starts with} $I$ if $I=J$ or $I\ll J$. Thus $I<J<I'$ with $I'$
starts with $I$ implies that $J$ starts with $I$.

Let $\Lambda$ be an artin algebra and  $\mod\Lambda$ the category of
finite generated (right) $\Lambda$-modules. For each
$M\in\mod\Lambda$, let $\mu(M)$ be the maximum of the sets
$\{|M_1|,|M_2|,\ldots, |M_t|\}$, where $M_1\subset M_2\subset \ldots
\subset M_t$ is a chain of indecomposable submodules of $M$. We call
$\mu(M)$ the {\bf Gabriel-Roiter measure}  of $M$. A subset $\mu$ of
$\mathcal{P}(\mathbb{N}_1)$ is called a GR measure if there is an
indecomposable module $M$ with $\mu(M)=\mu$. If $M$ is an
indecomposable $\Lambda$-module, we call an inclusion $X\subset M$
with $X$ indecomposable a {\bf GR inclusion} provided
$\mu(M)=\mu(X)\cup\{|M|\}$, thus if and only if every proper
submodule of $M$ has Gabriel-Roiter measure at most $\mu(X)$. In
this case, we call $X$ a {\bf GR submodule} of $M$.  Note that the
factor of a GR inclusion is indecomposable.

{\bf Remark.} We have seen in Introduction a different way to define
the Gabriel-Roiter measure. These two  definitions (orders) can be
identified. Namely, for each
$I=\{a_i|i\}\in\mathcal{P}(\mathbb{N}_1)$, let
$\mu(I)=\sum_{i}\frac{1}{2^{a_i}}$. Then $I<J$ if and only if
$\mu(I)<\mu(J)$.

\smallskip

We collect some results concerning Gabriel-Roiter measures, which
will be used later on.

\begin{prop}\label{GR} Let $\Lambda$ be an artin algebra and $X\subset M$ a GR
inclusion.
\begin{itemize}
  \item[(1)] If all irreducible maps to $M$ are monomorphisms, then
              the GR inclusion is an irreducible map.
  \item[(2)] Every irreducible map to $M/X$ is an epimorphism.
  \item[(3)] There is an irreducible monomorphism $X\ra Y$ with $Y$
             indecomposable and an epimorphism $Y\ra M$.
  \item[(4)] There is an epimorphism $\tau^{-1}X\ra M/X$.
\end{itemize}
\end{prop}
The proof of the above statements can be found for example in
\cite{Ch1,Ch2}.
\\
\\{\bf Example.} Let $Q$ be the Kronecker quiver:
$$\xymatrix{1\ar@/^/[r]\ar@/_/[r]& 2\\}.$$
We now describe the GR measures, which will be very useful in our
later discussion. The finite dimension representations (over an
algebraically closed field $k$) are simply called modules.

The GR measure of the indecomposable module with dimension vector
$(n,n+1)$ is $\{1,3,5,\ldots,2n+1\}$. The take-off modules are these
preprojective modules together with the simple injective module.

Every indecomposable regular module with dimension vector $(n,n)$
has GR measure $\{1,2,4,6,\ldots, 2n\}$. An indecomposable module is
in the central part if and only if it is a regular module.

The GR measure of the indecomposable module with dimension vector
$(n+1,n)$ is $\{1,2,4,\ldots, 2n,2n+1\}$. The landing part consists
of all non-simple indecomposable preinjective modules.

\subsection{A wild quiver}
From now on, we fix an algebraically closed field $k$ and consider
the following wild quiver $Q$:

$$\xymatrix{1\ar@{->}[r] & 2\ar@/^/[r]\ar@/_/[r]& 3\\}.$$

We refer to \cite{ARS,R3} for general facts of Auslander-Reiten
theory and refer to, for example, \cite{K,R0} for some basic results
of wild hereditary algebras.  The Cartan matrix $C$ and the Coxeter
transformation $\Phi$ are the following:
\begin{displaymath}
\begin{array}{ccc}
C=\left(\begin{array}{ccc}1 & 0 & 0\\
                          1 & 1 & 0\\
                          2 & 2 & 1\\ \end{array}\right), &
\Phi=-C^{-t}C=
\left(\begin{array}{rrr}0 & 1 & 0\\
                        3 & 3 & 2\\
                       -2 & -2 & -1\\\end{array}\right), &
\Phi^{-1}=
\left(\begin{array}{rrr}-1 & -1 & -2\\
                         1 & 0 & 0\\
                         0 & 2 & 3\\\end{array}\right).$$
\end{array}
\end{displaymath}

Then one many calculate the dimension vectors using  $\udim\tau
M=(\udim M)\Phi$ if $M$ is not projective and $\udim
\tau^{-1}N=(\udim N)\Phi^{-1}$ if $N$ is not injective, where $\tau$
denotes the Auslander-Reiten translation. The Eular form is
$\langle\underline{x},\underline{y}\rangle
=x_1y_1+x_2y_2+x_3y_3-x_1y_2-2x_2y_3$. Then for two indecomposable
modules $X$ and $Y$,
$$\dimHom(X,Y)-\dimExt^1(X,Y)=\langle\udim X,\udim Y\rangle$$

The Auslander-Reiten quiver consists of one preprojective component,
one preinjective component and infinitely many regular ones. An
indecomposable regular module $X$ is called quasi-simple if the
Auslander-Reiten sequence starting with $X$ has indecomposable
middle term. For each indecomposable regular module $M$, there is a
unique quasi-simple module $X$ and a unique natural number $r\geq 1$
(called quasi-length of $M$ and denoted by $\ql(M)=r$) such that
there is a sequence of irreducible monomorphisms $X=X[1]\ra
X[2]\ra\ldots\ra X[r]=M$.

The following is part of a regular component of the Auslander-Reiten
quiver containing an indecomposable module with dimension vector
$(0,1,1)$:

$$\xymatrix@C=0pt@R=8pt{
       {\tiny\left(\begin{array}{ccc}6 & 8 & 5\\\end{array}\right)}\ar@{->}[rd]
       &&{\tiny\left(\begin{array}{ccc}2 & 4 & 3\\\end{array}\right)}\ar@{->}[rd]
       &&{\tiny\left(\begin{array}{ccc}2 & 4 & 5\\\end{array}\right)}\ar@{->}[rd]
       &&{\tiny\left(\begin{array}{ccc}2 & 8 & 11\\\end{array}\right)}
       \\
       \cdots&{\tiny\left(\begin{array}{ccc}2 & 3 & 2\\\end{array}\right)}\ar@{->}[rd]\ar@{->}[ru]
       &&{\tiny\left(\begin{array}{ccc}1 & 2 & 2\\\end{array}\right)}\ar@{->}[rd]\ar@{->}[ru]
       &&{\tiny\left(\begin{array}{ccc}1 & 3 & 4\\\end{array}\right)}\ar@{->}[rd]\ar@{->}[ru]
       &\cdots\\
      {\tiny\left(\begin{array}{ccc}1 & 2 & 1\\\end{array}\right)}\ar@{->}[ru]
      &&{\tiny\left(\begin{array}{ccc}1 & 1 & 1\\\end{array}\right)}\ar@{->}[ru]
      &&{\tiny\left(\begin{array}{ccc}0 & 1 & 1\\\end{array}\right)}\ar@{->}[ru]
      &&{\tiny\left(\begin{array}{ccc}1 & 2 & 3\\\end{array}\right)}\\}$$

Let's denote by $H(1)$ an indecomposable module with dimension
vector $(0,1,1)$. Note that the indecomposable modules with
dimension vector $(0,1,1)$ are actually indexed by the projective
line $\mathbb{P}_k^1$.

\begin{lemm}\label{H}
{\rm (1)} For each $i\geq 0$, $\tau^{-i}H(1)$ contains no proper
regular submodules. In particular, a GR submodule of $\tau^{-i}H(1)$ is  preprojective. \\
{\rm (2)} For each $i\geq 0$, $\tau^{i}H(1)$ has no proper regular
factors. In particular, a GR factor module of $\tau^{i}H(1)$ is
preinjective.

\end{lemm}
\begin{proof} We show (1) and (2) follows similarly.
If $X$ is a proper regular submodule of  $\tau^{-i}H(1)$, then the
inclusion $X\subset \tau^{-i}H(1)$ induces a proper monomorphism
$\tau^{i}X\ra H(1)$. This is a contradiction since $H(1)$ has no
proper regular submodule.
\end{proof}

\begin{lemm}\label{moepi}
\begin{itemize}
   \item[(1)] There is a sequence of monomorphisms
              $$H(1)\ra \tau H(1)\ra\ldots\ra
           \tau^i H(1)\ra \tau^{i+1}H(1)\ra\ldots$$
   \item[(2)] There is a sequence of epimorphisms $$\ldots
             \tau^{-(i+1)}H(1)\ra \tau^{-i}H(1)\ra\ldots\ra \tau^{-1}H(1)\ra H(1).$$
\end{itemize}
\end{lemm}

\begin{proof}(1)
By above lemma,  a non-zero homomorphism from $\tau^i H(1)$, $i\geq
0$, to a regular module is a monomorphism. On the other hand,
\begin{displaymath}
  \begin{array}{cl}&\dimHom(H(1),\tau
                    H(1))-\dimExt^1(H(1),\tau H(1))\\
                  =& \langle\udim H(1), \udim \tau
                     H(1)\rangle\\
                  =&\langle(0,1,1),(1,1,1)\rangle\\
                  = &  0.\\
  \end{array}
\end{displaymath}
Since $\Ext^1(H(1), \tau H(1))\neq 0$, $\Hom(H(1),\tau H(1))\neq 0$
and thus there is monomorphism $H(1)\ra \tau H(1)$. It follows that
there is a sequence of monomorphisms
$$H(1)\ra \tau H(1)\ra\ldots\ra \tau^i H(1)\ra \tau^{i+1}H(1)\ra\ldots.$$

(2) follows similarly.
\end{proof}

It is easily seen that the inclusions $H(1)\ra\tau H(1)\ra\tau^2
H(1)$ are both GR inclusions. Thus the GR measures of $\tau^i H(1)$
are
$$\mu (H(1))=\{1,2\}, \mu(\tau H(1))=\{1,2,3\}, \mu(\tau^2H(1))=\{1,2,3,4\}$$
and $$\{1,2,3,4\}<\mu(\tau^iH(1))<\{1,2,3,4,5\}=\mu(I_3), \forall
i\geq 3.$$

\section{The Gabriel-Roiter submodules}

In this section, we will characterize the GR submodules of the
indecomposable preprojective modules, and the GR submodules  of the
indecomposable regular modules of the form $\tau^{-i}M,i\geq 0$  for
Kronecker modules $M$ (meaning that the indecomposable modules with
$(\udim M)_1=0$). We will also describe the take-off part and the
minimal central measure.

\subsection{The Gabriel-Roiter submodules of the preprojective modules}
We first calculate the GR submodules of the indecomposable
preprojetive modules. The dimension vectors of the indecomposable
projective modules are $\udim P_3=(0,0,1)$, $\udim P_2=(0,1,2)$ and
$\udim P_1=(1,1,2)$, and the beginning part of the preprojective
component is the following:
$$
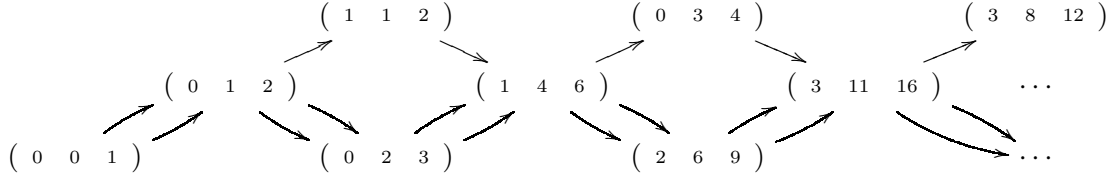
\xymatrix@C=0pt@R=8pt{
       &&{\tiny\left(\begin{array}{ccc}1 & 1 & 2\\\end{array}\right)}\ar@{->}[rd]
       &&{\tiny\left(\begin{array}{ccc}0 & 3 & 4\\\end{array}\right)}\ar@{->}[rd]
       &&{\tiny\left(\begin{array}{ccc}3 & 8 & 12\\\end{array}\right)}
       \\
       &{\tiny\left(\begin{array}{ccc}0 & 1 & 2\\\end{array}\right)}\ar@/^/[rd]\ar@/_/[rd]\ar@{->}[ru]
       &&{\tiny\left(\begin{array}{ccc}1 & 4 & 6\\\end{array}\right)}\ar@/^/[rd]\ar@/_/[rd]\ar@{->}[ru]
       &&{\tiny\left(\begin{array}{ccc}3 & 11 & 16\\\end{array}\right)}\ar@/^/[rd]\ar@/_/[rd]\ar@{->}[ru]
       &\cdots\\
      {\tiny\left(\begin{array}{ccc}0 & 0 & 1\\\end{array}\right)}\ar@/^/[ru]\ar@/_/[ru]
      &&{\tiny\left(\begin{array}{ccc}0 & 2 & 3\\\end{array}\right)}\ar@/^/[ru]\ar@/_/[ru]
      &&{\tiny\left(\begin{array}{ccc}2 & 6 & 9\\\end{array}\right)}\ar@/^/[ru]\ar@/_/[ru]
      &&\cdots\\}$$

\begin{prop}\label{P}
  \begin{itemize}
     \item[(1)] Up to isomorphism, $\tau^{-i}P_2$ is the unique Gabriel-Roiter submodule of
                $\tau^{-(i+1)}P_3$ for each $i\geq 0$.
     \item[(2)] Up to isomorphism, $\tau^{-i}P_3$ is the unique Gabriel-Roiter submodule of
                $\tau^{-i}P_1$ for each $i\geq 1$.
     \item[(3)] A Gabriel-Roiter submodule of
                $\tau^{-i}P_2\cong\left\{\begin{array}{ll}\tau^{-(i-1)}P_1, &  \textrm{$i$ is odd};\\
                                               \tau^{-i}P_3,    &  \textrm{$i$ is even}.\\
                            \end{array}
                       \right.$
   \end{itemize}
\end{prop}

\begin{proof}
(1) Since $\tau^{-i}P_2\oplus \tau^{-i}P_2\ra  \tau^{-(i+1)}P_3\ra
0$ is a right minimal almost split morphism and the irreducible maps
involved are monomorphisms, a GR submodule of $\tau^{-(i+1)}P_3$ is
isomorphic to $\tau^{-i}P_2$ (Proposition \ref{GR}(1)).

(2) We first show that $\mu(\tau^{-r}P_3)>\mu(X)$ for all
predecessors $X$ in the preprojective component.  This is obvious
for $r=1$. Since the irreducible maps
$\tau^{-(r-1)}P_1\ra\tau^{-r}P_2\ra\tau^{-(r+1)}P_3$ are all
monomorphisms,
$\mu(\tau^{-(r+1)}P_3)>\mu(\tau^{-r}P_2)>\mu(\tau^{-(r-1)}P_1)$.
Note that a predecessor of  $\tau^{-(r+1)}P_3$ is either isomorphic
to $\tau^{-r}P_2$ or  $\tau^{-(r-1)}P_1$, or is a predecessor of
$\tau^{-r}P_3$.  Therefore, we can finish the proof by induction.

Now we show that $\tau^{-i}P_3$ is a GR submodule of $\tau^{-i}P_1$
for each $i\geq 1$. Since there is a sectional path $\tau^{-i}P_3\ra
\tau^{-i}P_2\ra\tau^{-i}P_1$, thus the composition of irreducible
maps  is either an epimorphism or a monomorphism. But for each
$i\geq 1$,
\begin{eqnarray}
    |\tau^{-i}P_3|-|\tau^{-i}P_1| &=& |\tau^{-i}P_3|+|\tau^{-(i-1)}P_1|-|\tau^{-i}P_2|\nonumber\\
         &=&  (2|\tau^{-i}P_3|+|\tau^{-(i-1)}P_1|-|\tau^{-i}P_2|)-|\tau^{-i}P_3| \nonumber \\
         &=&  |\tau^{-(i-1)}P_2|-|\tau^{-i}P_3|<0.\nonumber
\end{eqnarray}
Thus there is a monomorphism from  $\tau^{-i}P_3$ to $\tau^{-i}P_1$.
Since a GR submodule of $\tau^{-i}P_1$ is one of it's predecessors,
it is sufficient to show that neither $\tau^{-i}P_2$ nor
$\tau^{-(i-1)}P_1$ is a GR submodule of $\tau^{-i}P_1$. But this is
obvious since the irreducible map $\tau^{-i}P_2\ra \tau^{-i}P_1$ is
an epimorphism and $\Hom(\tau^{-(i-1)}P_1,\tau^{-i}P_1)=0$.

(3) Since all irreducible maps to $\tau^{-i}P_2$ are monomorphisms,
a GR submodule of $\tau^{-i}P_2$ is isomorphic to either
$\tau^{-(i-1)}P_1$ or $\tau^{-i}P_3$. First, we show that for each
$i\geq 1$, $\tau^{-(i-1)}P_1\ra\tau^{-i}P_2$ is a GR inclusion
implies that $\tau^{-(i+1)}P_3\ra \tau^{-(i+1)}P_2$ is a GR
inclusion.  If this is not the case, then $\tau^{-i}P_1\ra
\tau^{-(i+1)}P_2$ is a GR submodule. Then we have
$$\mu(\tau^{-i}P_3)<\mu(\tau^{(-i-1)}P_1)<\mu(\tau^{-i}P_2)<\mu(\tau^{-(i+1)}P_3)<
   \mu(\tau^{-i}P_1)<\mu(\tau^{-(i+1)}P_2).$$
Since $\tau^{-i}P_3$ is a GR submodule of $\tau^{-i}P_1$,
$\mu(\tau^{-(i-1)}P_1)$ starts with $\mu(\tau^{-i}P_3)$. In
particular, there exists a submodule $X$ of $\tau^{-(i-1)}P_1$ such
that $\mu(X)=\mu(\tau^{-i}P_3)$. Since $X$ is not isomorphic to
$\tau^{-i}P_3$, $X$ has to be a predecessor of $\tau^{-i}P_3$ and
thus $\mu(X)<\mu(\tau^{-i}P_3)$ by the discussion in (2). This is a
contradiction.

Secondly, we show that   $\tau^{-i}P_3\ra\tau^{-i}P_2$ is a GR
inclusion implies that $\tau^{-i}P_1\ra \tau^{-(i+1)}P_2$ is a GR
inclusion.  It is sufficient to show that
$\tau^{-(i+1)}P_3<\tau^{-i}P_1$. We may assume $i\geq 1$. Then
$\tau^{-i}P_3$ is also a GR submodule of $\tau^{-i}P_1$. Since the
irreducible map $\tau^{-i}P_2\ra \tau^{-i}P_1$ is an epimorphism,
$\mu(\tau^{-i}P_1)>\tau^{-i}P_2$. Note that $\tau^{-i}P_2$ is a GR
submodule of $\tau^{-(i+1)}P_3$. It follows that
$\mu(\tau^{-(i+1)}P_3)<\mu(\tau^{-i}P_1)$.

Now the statement follows by induction and the facts that $P_3$ is a
GR submodule of $P_2$ and $P_1$ is a GR submodule of $\tau^{-1}P_2$.
\end{proof}

The following observations can be easily checked:
\begin{itemize}
      \item  $\mu(\tau^{-i}P_2)>\mu(X)$ if $X$ is a predecessor of  $\tau^{-i}P_2$ for all $i\geq 0$.
      \item  $\mu(\tau^{-i}P_1)>\mu(X)$ for all predecessor $X$ if
$\left\{\begin{array}{ll}\textrm{ $i$ is even};\\
                      \textrm{ $i$ is odd and}\,\ X\ncong \tau^{-(i-1)}P_1,\tau^{-i}P_2.
\end{array}\right.$
\end{itemize}

\subsection{The take-off part and the minimal central measure}

Let $P_1$ be the indecomposable projective module, i.e., $\udim
P_1=(1,1,2)$. If $X$ is a non-injective indecomposable proper factor
of $P_1$, then $X$ has dimension vector  $\udim X=(1,1,1)$ and
$\mu(X)=\{1,2,3\}$. Thus a non-simple indecomposable module $M\ncong
P_1$ with $\Hom(P_1,M)\neq 0$ has GR measure
$\mu(M)>\mu(P_1)=\{1,3,4\}$.

\begin{lemm} Let $M$ be an indecomposable module, which is neither simple nor
injective.
   \begin{itemize}
       \item[(1)] The GR measure of $M$ is
                  $\{1,3,5,\ldots,2n+1\}$ if and only if $\udim
                  M=(0,n,n+1)$.
       \item[(2)] The GR measure of $M$ is
                  $\{1,2,4,\ldots,2n\}$ if and only if $\udim M=(0,n,n)$.
   \end{itemize}
\end{lemm}

\begin{proof}We show (1) and (2) follows similarly. By the
description of the GR measures of Kronecker modules, $\udim
M=(0,n,n+1)$ implies that $\mu(M)=\{1,3,5,\ldots,2n+1\}$. For the
converse implication, we use induction on the length. It is clear
that $\mu(M)=\{1,3\}$ if and only if $M$ is the projective module
$P_2$, i.e., $\udim M=(0,1,2)$. Now assume that
$\mu(M)=\{1,3,5,\ldots,2n+1,2n+3\}$ with $n\geq 1$. Then a GR
submodule $X$ of $M$ has GR measure $\{1,3,5,\ldots,2n+1\}$. Thus by
induction $\udim X=(0,n,n+1)$. Since the GR factor module $M/X$ has
length $2$, its dimension vector $\udim M/X=(1,1,0)$ or $(0,1,1)$.
In the first case, $M/X$ is the indecomposable injective module
$I_2$. However, $I_2$ cannot be a GR factor module, since there is
an irreducible monomorphism $S_2\ra I_2$ (Proposition \ref{GR} (2)).
It follows that $\udim M/X=(0,1,1)$ and $\udim M=(0,n+1,n+2)$.
\end{proof}

\begin{prop}\label{takeoff} A non-simple indecomposable module $M$ is a take-off
module if and only if $\udim M=(0,n,n+1)$ for some $n\geq 1$. Thus
the take-off measures are of the form $\{1,3,5,\ldots,2n+1\}$ for
$n\geq 0$.
\end{prop}
\begin{proof} Let $\mu_n=\{1,3,5,\ldots,2n+1\}$.
It is sufficient to show that $\mu_{n+1}$ is a direct successor of
$\mu_n$ for each $n\geq 0$. Assume for a contradiction that
$$\{1,3,\ldots,2n+1\}=\mu_n<\mu<\mu_{n+1}=\{1,3,\ldots,2n+1,2n+3\}.$$
Then $\mu=\{1,3,\ldots,2n+1,m_1,m_2,\ldots,m_s\}$ with $m_1>2n+3$.
In particular, there exists an indecomposable module $X$, containing
$Y$ with $\udim Y=(0,n,n+1)$ as a GR submodule and the corresponding
GR factor $X/Y$ has length greater than $2$. Assume that $\udim
X/Y=(a,b,c)$. By the description of the GR measures of the Kronecker
modules, we have $a\neq 0$ and thus $\Hom(P_1,X)\neq 0$. Therefore
either there is a monomorphism $P_1\ra X$, or $X$ contains an
indecomposable submodule with dimension vector $(1,1,1)$. It follows
that $\mu\geq\mu(X)\geq \mu(P_1)>\mu_r$ for all $r\geq 0$. This
contradiction implies that $\mu_{n+1}$ is a direct successor of
$\mu_n$ for each $n\geq 1$.
\end{proof}

\begin{prop}\label{minimal} The indecomposable projective module $P_1$ is a central
module and $\mu(P_1)$ is the minimal central measure. In particular,
$\mu(P_1)$ has no direct predecessor.
\end{prop}

\begin{proof}Since $\mu(P_1)=\{1,3,4\}>\mu_r=\{1,3,5,\ldots, 2r+1\}$ for all
$r\geq 0$, $P_1$ is a central module. Assume for a contradiction
that $\mu$ is a central measure such that $\mu<\mu(P_1)$. Then
$$\{1,3,5,\ldots,2r+1\}=\mu_r<\mu<\mu(P_1)$$ since $\mu_r$ is a take-off measure
for each $r\geq 0$. It follows that $\mu$ starts with
$\{1,3,5,\ldots,2n+1,2n+2\}$ for some $n\geq 2$. Therefore, there is
an indecomposable module $M$ with length $2n+2$ and containing the
indecomposable module $(0,n,n+1)$ as a GR submodule. If the
dimension vector of $M$ is $(1,n,n+1)$, then $\Hom(P_1,M)\neq 0$ and
$\mu(M)>\mu(P_1)$. Thus the only possibility is that $\udim
M=(0,n+1,n+1)$ and the GR measure is $\{1,2,4,\ldots,2n\}$ by
previous proposition. This is  a contradiction. Therefore,
$\mu(P_1)$ is the minimal central measure and thus has no direct
predecessor.
\end{proof}

\subsection{The Gabriel-Roiter submodules of $\tau^{-i}M$ with $M$ a Kronecker module}

For each integer $a>0$ and $\lambda\in\mathbb{P}^1_k$, we denote by
$H(a)_\lambda$ the indecomposable module with dimension vector
$(0,a,a)$ and parameter $\lambda$, and  by $H_a$ and $H^a$ for
$a\geq 0$, the unique indecomposable module with dimension vector
$(0,a,a+1)$ and $(0,a+1,a)$, respectively. An indecomposable module
$X$ is a submodule (resp. factor module) of $H_a$ (resp. $H^a$) if
and only if $X$ is isomorphic to $H_b$ (resp. $H^b$) for some $b\leq
a$.  Similarly, an indecomposable module $Y$ is a submodule (resp.
factor module) of $H(a)_\lambda$  if and only if $Y$ is isomorphic
to $H_b$ or $H(b)_\lambda$ (resp. $H^b$ or $H(b)_\lambda$) for some
$b\leq a$.   Note that the GR measures are
$\mu(H_a)=\{1,3,5,7,\cdots,2a+1\}$,
$\mu(H(a)_\lambda)=\{1,2,4,6\cdots, 2a\}$ and
$\mu(H^a)=\{1,2,4,6,\cdots,2a,2a+1\}$.

\begin{lemm}
\begin{itemize}
   \item[(1)] $H(a)_\lambda$ is a quasi-simple module for each $a\geq 1$ and $\lambda\in\mathbb{P}_k^1$.
   \item[(2)] $H_a$ is a quasi-simple module for each $a\geq 4$ and $H^a$
                is
             quasi-simple module for each $a\geq 1$.
   \item[(3)] Any two regular
             modules of above three kinds are in different regular components except the
           pair $(H^1,H_4)$, where $\tau^2H_4=H^1$.
\end{itemize}
\end{lemm}

\begin{proof} We show (1) and (2) follows similarly. Assume  that
$H(a)_\lambda$ is not a quasi-simple module. Then there is a
quasi-simple module $X$ and an integer $r\geq 2$ such that
$X[r]=H(a)_\lambda$. Then $X\cong H_b$ or $H(b)_\lambda$ for some
$0<b<a$. Thus $\udim X=(0,b,b+1)$ or $\udim X=(0,b,b)$. However, the
dimension vector of $\tau^{-1}X$ is
     $$(0,b,b')\left(\begin{array}{rrr}-1 & -1 & -2\\
                                       1 & 0  & 0\\
                                       0 & 2  & 3
               \end{array}\right)=(b,2b',3b')$$
for $b'=b$ or $b+1$. It follows that $(\udim H(1)_\lambda)_1=(\udim
X[r])_1\geq b$. This is a contradiction. Therefore $H(a)_\lambda$ is
quasi-simple.

Now we prove (3).   It is clear that $(\udim \tau
H(1)_\lambda)_1\neq 0\neq (\udim \tau^{-1}H(1)_\lambda)$. It follows
from Lemma \ref{moepi} that $(\udim \tau^iH(1)_\lambda)_1\neq 0$ for
any $i\neq 0$. In particular, $\tau^iH(1)_\lambda$ does not
isomorphic to $H_a$ or $H^a$ for any $a$, or $H(b)_\gamma$ for any
$b>1$ or $b=1$ and $\gamma \neq \lambda$. If $a\geq 2$, the short
exact sequence $0\ra H(1)_\lambda\ra H(a)_\lambda\ra H(a-1)_\lambda
\ra 0$ induces an exact sequence  $0 \ra \tau^iH(1)_\lambda\ra
\tau^i H(a)_\lambda\ra \tau^i H(a-1)_\lambda\ra 0$ for each integer
$i$. Thus ($\udim
\tau^iH(a)_\lambda)_1>(\udim\tau^iH(1)_\lambda)_1\geq 1$ for every
$i\geq 0$. It follows that $H(a)_\lambda$ is the unique
indecomposable module such that $(\udim M)_1=0$ in the component
containing $H(a)_\lambda$.

In stead of $H(a)_\lambda$, we simply write $H(a)$ in the following
proof,  since the parameter is not so important. We show that $H^a$
and $H^b$ (similarly $H_a$ and $H_b$) are not in the same component.
Without loss of generality, we may assume that $\tau^iH^b=H^a$ for
some $i>0$. Since $H(b)$ is a submodule $H^b$, $\tau^iH(b)$ is thus
a submodule of $\tau^iH^b=H^a$. Thus $\tau^iH(b)$ is isomorphic to
$H(c)$ or $H_c$ for some $c<a$. It follows that this $H(b)$ and
$H(c)$ or $H_c$ are in the same component.  This is a contradiction.

To finish the proof, it is sufficient to show that  $H^a$ and $H_b$
are not in the same component with only one exception. If
$H_b=\tau^iH^a$ for some $i>0$, then as before, we have a
monomorphism $\tau^iH(a)\ra \tau^iH^a=H_b$. Thus $\tau^iH(a)$ is
isomorphic to $H_c$ for some $c<b$. This is a contradiction since
$H(a)$ and $H_c$ are not in the same component. Therefore the only
possibility is $H^a=\tau^iH_b$ for some $i>0$. If $b>4$, $H_{b-1}$
is a regular submodule of $H_b$ with factor module $H(1)$. Thus
$\tau^iH_{b-1}$ is a submodule of $\tau^iH_b=H^a$ with factor
$\tau^iH(1)$. Since any indecomposable factor module of $H^a$ is of
the form $H^c$ with $c<a$,  $\tau^iH(1)$ is isomorphic to some
$H^c$. This is again a contradiction. Thus $H^a=\tau^iH_b$ may
happen only in case $b=4$. An easy calculation shows that
$\tau^2H_4=H^1$.
\end{proof}

\begin{lemm}\label{ql} Let $M$ be an indecomposable regular module with
dimension vector $\udim M=(a,b,c)$. Then the quasi-length of $M$
satisfies $ql(M)\leq a+1$. Moreover, if $a=1$ and $ql(M)=2$, then
$\udim M=(1,2,2)$ or $(1,3,4)=(1,2,2)\Phi^{-1}$.
\end{lemm}

\begin{proof}
Assume for a contradiction that $M=X[r]$ for some quasi-simple
module $X$ and $r\geq a+2$. Then
$\sum_{i=0}^{r-1}(\udim\tau^{-i}X)_1=a\leq r-2$. It follows from
above discussion that there are  $0\leq i<j\leq r-1$ such that
$(\udim(\tau^{-i}X))_1=0=(\udim(\tau^{-j}X))_1$ and
$(\udim(\tau^{-s}X))_1=1$ for all $0\leq s\neq i,j\leq r-1$. The
only possibility is that $\udim(\tau^{-i}X)=(0,2,1)$ and
$\udim(\tau^{-j}X)=(\udim(\tau^{-(i+2)}X)=(0,4,5)$. But in this case
$(\udim(\tau^{-(i+1)}X)_1=(2,2,3)$, which contradicts
$(\udim(\tau^{-(i+1)}X))_1=1$.

If $a=1$ and $\ql(M)=2$, then $(\udim X)_1=0$ or $(\udim
\tau^{-1}X)_1=0$. If $\udim X=(0,a,b)$, then $\udim
\tau^{-1}X=(a,2b,3b)$. It follows that $a=1$, $b=1$ and thus $\udim
M=(1,3,4)$.  If $\udim \tau^{-1}X=(0,a,b)$, then $\udim
X=(3a-2b,3a-2b,2a-b)$.  Thus $3a-2b=1$ and  only possibility is
$a=b=1$. It follows that $\udim M=(1,2,2)$.  Note that
$(1,2,2)=(1,3,4)\Phi$.
\end{proof}

\begin{prop}\label{H_a} Up to isomorphism,  $\tau^{-i}H_a$ is the unique
GR submodule of $\tau^{-i}H_{a+1}$ for each $a\geq 1$ and $i\geq 0$.
It follows that all $\tau^{-i}H(1)_\lambda$ are GR factor modules.
\end{prop}

\begin{proof}
Since $H_b$ is a GR submodule of $H_{b+1}$ with GR factor module
$H(1)_\lambda$. (Different embeddings give rise to different
factors.) Thus we have monomorphisms  $\tau^{-i}H_b\ra
\tau^{-i}H_{b+1}$ with factors $\tau^{-i}H(1)_\lambda$.  In
particular,
$\mu(\tau^{-i}H_b)<\mu(\tau^{-i}H_{b+1})<\ldots<\mu(\tau^{-i}H_a)$
for all $b<a$.

Assume that $X$ is  a GR submodule of $\tau^{-i}H_{a+1}$. If $X$ is
a regular module, then the monomorphism $X\ra \tau^{-i}H_{a+1}$
induces a monomorphism $\tau^{i}X\ra H_{a+1}$. It follows that
$\tau^{i}X\cong H_b$ for some $b\leq a$. However,
$\mu(\tau^{-i}H_b)<\mu(\tau^{-i}H_{b+1})$ for all $b$. Thus
$X\cong\tau^{-i}H_a$. Similarly, if $X\cong \tau^{-j}P$ is
preprojective for some indecomposable projective module $P$ and some
$j>i+1$, then $\tau^{-(j-i)}P\cong H_b$ and thus $j-i\leq 1$, which
is impossible. Thus $X\cong \tau^{-j}P$ for some indecomposable
module $P$ and some $j\leq i+1$.

For $a=1,2$, all modules  $\tau^{-i}H_{a+1}$ are preprojective and
the statement is the same with that $\tau^{-i}P_2$ is a GR submodule
of $\tau^{-(i+1)}P_3$ and  $\tau^{-i}P_3$ is a GR submodule of
$\tau^{-i}P_1$, which we have proved (Proposition \ref{P}). If
$a=3$, then the GR submodule $X$ of $\tau^{-i}H_4$ has to be
preprojective since $H_4$ contains no regular submodules. In this
case,  $X\cong \tau^{-j}P$ for some projective module $P$ and some
$j\leq i+1$. Note that $X$ is a predecessor of $\tau^{-i}H_3$ in the
preprojective component. If $i$ is odd, then
$\mu(\tau^{-i}H_3)=\mu(\tau^{-(i+1)}P_1)$ is larger than all
$\mu(Y)$ if $Y$ is one of its predecessors.  Thus $\tau^{-i}H_3\ra
\tau^{-i}H_4$ is a GR inclusion. If $i$ is even, then the only
predecessors of $\tau^{-i}H_3$ with  GR measure larger than
$\mu(\tau^{-i}H_3)$ are $\tau^{-(i+1)}P_2$ and $\tau^{-i}P_1$. In
both cases, we get  monomorphisms from $\tau^{-1}P_2$ and $P_1$ to
$H_4$, respectively. This is a contradiction. Therefore,
$\tau^{-i}H_3$ is a GR submodule of $\tau^{-i}H_4$.

Finally, assume that $a\geq 4$. It is sufficient to show that the GR
submodule of $\tau^{-i}H_a$ is regular for each $i\geq 0$. If $X$ is
preprojective, as before, $X$ is a predecessor of $\tau^{-i}H_3$.
Again if $i$ is odd, then $\mu(X)\leq
\mu(\tau^{-i}H_3)<\mu(\tau^{-i}H_a)$, a contradiction. If $i$ is
even, we may repeating the arguments as in $a=3$ case and get a
contradiction.

We finish the proof.
\end{proof}

{\bf Remark.}  Since $\tau^{-i}H(1)_\lambda$ are GR factors, we may
obtain, for each natural number $r$, that a GR inclusion $X\ra Y$
such that $|Y/X|\geq r$. However, for tame quivers, the dimension
vectors of the GR factors are always bounded by $\delta$, where
$\delta$ is the minimal positive imaginary root \cite{Ch2}.

\begin{prop}\label{H(1)} Fix a $\lambda\in\mathbb{P}_k^1$ and we simply denote $H(1)_\lambda$ by $H(1)$.
For each $i>0$, the GR submodule of
$\tau^{-i}H(1)=\left\{\begin{array}{ll}\tau^{-(i-1)}P_1, &  \textrm{$i$ is odd};\\
                                      \tau^{-(i-1)}P_2, & \textrm{$i$ is even}.\end{array}\right.$
\end{prop}

\begin{proof} Let $X$ be a GR submodule of $\tau^{-i}H(1)$.
Then $X$ is preprojective by Lemma \ref{H}.  If $X=\tau^{-j}P$ for
some indecomposable projective module $P$ and $j>i$,  then  we
obtain a monomorphism from $\tau^iX$ to $H(1)$. But this is
impossible since the unique proper submodule of $H(1)$ is a the
simple projective module. Thus $X=\tau^{-r}P$ for some
indecomposable projective module $P$ and some $r<i$.  Clearly, $P_1$
is a GR submodule of $\tau^{-1}H(1)$ with $H(1)$ as GR factor. We
thus obtain monomorphisms $\tau^{-(i-1)}P_1\ra \tau^{i}H(1)$ for all
$i>0$.  By the same reason, we get monomorphisms $\tau^{-i}P_2\ra
\tau^{i}H(1)$. We can finish the proof by applying the description
of the GR measures of the preprojective modules.
\end{proof}

\begin{coro}Fix a $\lambda\in\mathbb{P}_k^1$.
Then $\mu(P)<\mu(\tau^iH(1))<\mu(\tau^jH(1))$ for all $i<j$ and all
preprojective module $P$.
\end{coro}
\begin{proof}
We need only check that $\mu(\tau^{-i}H(1))>\mu(\tau^{-(i+1)}H(1))$
for all $i\geq 0$. This is clear for $i=0$.  Now assume $i>0$ is
odd. We use the following diagram to indicate the homomorphisms:
$$\xymatrix{
 \tau^{-(i-1)}P_1\ar@{->}[r]^{\textrm{GR}}\ar@{->}[d]_{\textrm{GR}} & \tau^{-i}H(1)\\
 \tau^{-i}P_2\ar@{->}[r]_{\textrm{GR}}\ar@{->}[ru]_{\textrm{epi.}}
 & \tau^{-(i+1)}H(1)\ar@{->}[u]_{\textrm{epi.}}\\}$$
Thus $\mu(\tau^{-i}H(1))>\mu(\tau^{-i}P_2)>\mu(\tau^{-(i+1)}H(1))$.
The case that $i>0$ is even follows similarly. Then using the
description of the GR submodules of $\tau^{-i}H(1)_\lambda$ and
those of preprojective modules, we may easily deduce that
$\mu(\tau^{-i}H(1))>\mu(P)$ for all $i$ and all preprojective
modules $P$.
\end{proof}

\begin{prop}\label{H(a)}Fix a $\lambda\in\mathbb{P}_k^1$. For each $i\geq 0$ and $a\geq 1$,
$\tau^{-i}H(a)$ is the unique, up to isomorphism, GR submodule of
$\tau^{-i}H(a+1)$.
\end{prop}
\begin{proof} Since there is a monomorphism $\tau^{-i}H(1)\ra \tau^{-i}H(a+1)$,
the above corollary implies that the GR submodules of
$\tau^{-i}H(a)$ are regular modules for all $i\geq 0$. Let $X$ be a
GR submodule of $\tau^{-i}H(a)$. Then $\tau^{i}X$ is a submodule of
$H(a+1)$ and thus isomorphic to $H(b)$ or $H_b$ for some $b<a+1$.
Thus $X\cong\tau^{-i}H(b)$ or $X\cong \tau^{-i}H_b$.

Since there is a monomorphism from $\tau^{-i}H(a)$ to
$\tau^{-i}H(a+1)$, it is sufficient to show that
$\mu(\tau^{-i}H(a))\geq \mu(X)$. This is obvious for
$X\cong\tau^{-i}H(b)$. Assume that $X\cong\tau^{-i}H_b$ for some
$b<a+1$. We consider the dimension vectors of the form $(0,c,c+1)$.
Note that
$(0,c,c+1)\Phi^{-1}=(0,c,c)\Phi^{-1}+(0,0,1)\Phi^{-1}=(0,c,c)\Phi^{-1}+(0,2,3)$.
This implies $|\tau^{-1}H_c|> |\tau^{-1}H(c)|$. Since $(0,2,3)=\udim
H_2$, $(0,2,3)\Phi^{-i}$ is a positive vector, namely $\udim
\tau^{-i}H_2$ for each $i\geq 0$.  It follows that
$|\tau^{-i}H_c|>|\tau^{-i}H(c)|\geq|\tau^{-i}H(1)|$. Note that
$\mu(\tau^{-i}H_3)<\mu(\tau^{-i}H(1))$. Since $\tau^{-i}H_c$ is a GR
submodule of $\tau^{-i}H_{c+1}$ and
$|\tau^{-i}H_c|>|\tau^{-i}H(c)|>|\tau^{-i}H(1)|$ for all $c\geq 1$,
we have  $\mu(\tau^{-i}H_c)<\mu(\tau^{-i}H(1))\leq
\mu(\tau^{-i}H(a))$. Thus for $a\geq 1$, $\tau^{-i}H(a)$ is a GR
submodule of $\tau^{-i}H(a+1)$.
\end{proof}

\begin{prop}\label{H^a} For each $a\geq 1$, the GR submodules of
$\tau^{-i}H^a$, $i\geq 0$
are $$\left\{\begin{array}{ll}\tau^{-i}H(a)_\lambda,\lambda\in\mathbb{P}^1_k,  & i=0,1;\\
                               \tau^{-(i-1)}P_1, & i\geq 2,a=1;\\
                              \tau^{-i}H(a-1)_\lambda,\lambda\in\mathbb{P}^1_k,
                                 & i\geq 2,a\geq 2.\\ \end{array}\right.$$
\end{prop}

\begin{proof} First assume that $a=1$. Note that $\tau^{-2} H^1$
is $H_4$ and  the GR submodules of $\tau^{-i}H^1$ is already known
as $\tau^{-(i-1)}P_1$ for all $i\geq 2$. Every indecomposable module
$H(1)_\lambda$  is a GR submodule of $H^1$ with factor $S_2$, which
is not injective. It follows that there is a monomorphism
$\tau^{-1}H(1)_\lambda\ra \tau^{-1}(H^1)$. We claim that it is
namely a GR inclusion. Note that $\tau^{-1}H(1)_\lambda$ has GR
measure $\{1,3,4,6\}>\mu(P)$ for any preprojective module $P$. Thus
a GR submodule of $\tau^{-1}H^1$ has to be a regular one. Assume $X$
is a GR submodule of $\tau^{-1}(H^1)$. We obtain a monomorphism
$\tau X\ra H^1$, and therefore, $X$ has dimension $(0,1,1)$.  Thus
$\tau^{-1}H(1)_\lambda$ is a GR submodule of $\tau^{-1}H^1$ for each
$\lambda$.

Now we consider the case $a\geq 2$. Since there is a monomorphism
$H(1)_\lambda\ra H^a$ for each $\lambda$ with indecomposable regular
factor $H^{a-1}$, the GR submodules of $\tau^{-i}H^a$ are regular.
Let $X$ be a GR submodule of $\tau^{-i}H^a$, then there is a
monomorphism $\tau^iX\ra H^a$. It follows $X$ is either isomorphic
to $\tau^{-i}H_b$ or $\tau^{-i}H(b)_\lambda$ for some $b\leq a$.
From the description of the GR submodules  of these modules, we know
that the GR submodules of $\tau^{-i}H^a$ are of the form
$\tau^{-i}H(b)_\lambda$ with $b$ as large as possible. We may
calculate the dimension vectors as follows:
\begin{displaymath}
\begin{array}{rcl}
(0,a+1,a)\Phi^{-2}   & = & (a+1,2a,3a)\Phi^{-1}=(a-1,5a-1,7a-2)\\
(0,a,a)\Phi^{-2}     & = & (7a,2a,3a)\Phi^{-1}=(a,5a,7a)\\
(0,a-1,a-1)\Phi^{-2} & = & (a-1,2a-2,3a-3)\Phi^{-1}=(a-1,5a-5,7a-7)
\end{array}
\end{displaymath}
Comparing the dimension vectors, we conclude that the GR submodules
of $\tau^{-i}H^a$ for $a\geq 2$ and $i\geq 2$ are
$\tau^{-i}H(a-1)_\lambda$ for all $\lambda\in\mathbb{P}^1_k$, and
the GR submodules of $\tau^{-i}H^a$ are $\tau^{-i}H(a)_\lambda$ for
$i=0,1$.
\end{proof}

In general, for each $i\geq 0$, $\tau^iH(a)_\lambda$ (resp.
$\tau^iH^a$, $\tau^iH_a$) is not a GR submodule of
$\tau^iH(a+1)_\lambda$ (resp. $\tau^iH^{a+1}$, $\tau^iH_{a+1}$).
Next we will give one example.

\begin{prop}
Any GR factor of $\tau H(a)_\lambda$ is isomorphic to the simple
injective module. In particular, $\tau H(a)_\lambda$ is not a GR
submodule of $\tau H(a+1)_\lambda$.
\end{prop}
\begin{proof} This is obvious for $a=1$.  Now assume that $a>1$.
As before, it is easily seen that  a GR submodule of $\tau^i
H(a)_\lambda$ is regular  and has GR measure not smaller than
$\mu(H(1)_\lambda)$.

Let $0\ra X\ra \tau H(a)_\lambda\ra Y\ra 0$ be a GR sequence. If $Y$
is not injective, we get the following exact sequence $0\ra\tau^{-1}
X\ra H(a)_\lambda\ra \tau^{-1}Y\ra 0$. Then $\tau^{-1}X$ is
isomorphic to $H(b)$ or $H_b$ for some $b<a$. Because  $H_b$ is
cogenerated by $H(1)_\lambda$, $\tau H_b$ is cogenerated by $\tau
H(1)_\lambda$. Thus $\mu(\tau H_b)<\mu(\tau(H(1)_\lambda)$. Since
there is a monomorphism $\tau H(1)_\lambda\ra \tau H(a)_\lambda$,
$\tau^{-1}X$ is not of the form $H_b$. Assume that
$\tau^{-1}X=H(b)_\lambda$ for some $b<a$ and therefore $b=a-1$ since
$X$ is a GR submodule.

However, an easily calculation shows  $\udim \tau
H(a)_\lambda=(a,a,a)$ and
\begin{eqnarray}
\udim (\tau H(a-1)_\lambda)[2] & = & \udim\tau^{-1}H(a-1)_\lambda+\udim H(a-1)_\lambda\nonumber\\
                         & =& (0,a-1,a-1)+(a-1,a-1,a-1)\nonumber\\
                        & = & (a-1,2a-2,2a-2)\nonumber
\end{eqnarray}
Thus, there does not exist an epimorphism from $(\tau
H(a-1)_\lambda)[2]$ to $\tau H(a)_\lambda$ (Proposition
\ref{moepi}(3)). This contradiction implies that in the GR sequence
concerning $\tau H(a)_\lambda$, $Y$ has to be injective. It follows
that $Y$ is isomorphic to $I_1$ or $I_3$.

If $Y$ is isomorphic to $I_1$, then $\udim
X=(a,a,a)-(2,2,1)=(a-2,a-2,a-1)$ and then $\udim
\tau^{-1}X=(0,a,a+1)$. This is impossible since $(\udim
X[2])_1=a-2<a=(\udim \tau H(a)_\lambda)_1$. Thus the GR factor of
$\tau H(a)$ is the simple injective module $I_1$ and  the GR
submodule of $\tau H(a)$ has dimension vector $(a-1,a,a)$. In
particular, $\tau H(a)_\lambda$ is not a GR submodule of $\tau
H(a+1)_\lambda$ because $\udim \tau H(a)_\lambda=(a,a,a)\neq
(a,a+1,a+1)$.
\end{proof}

\subsection{The Gabriel-Roiter measures of indecomposable modules of small
dimensions}

\begin{lemm}
   \begin{itemize}
      \item[(1)] An indecomposable module $M$ has GR measure
                 $\{1,2,3\}$ if and only $\udim M=(1,1,1)$ or
                 $(0,2,1 )$.
      \item[(2)] An indecomposable module $M$ has GR measure
                 $\{1,2,3,4\}$ if and only if $\udim M=(1,2,1)$.
      \item[(3)] An indecomposable module $X$ with dimension vector
                $(1,2,2)$ and $\ql(M)=2$ has GR measure
                $\{1,2,3,5\}$.
      \item[(4)] An indecomposable module $M$ with dimension vector
                $(1,3,4)$ and $\ql(M)=2$ has GR measure
                $\{1,2,8\}$.
      \item[(5)] An indecomposable module $M$ has GR measure
             $\{1,3,4,6\}$ if and only if $\udim M=(1,2,3)$, i.e., $M\cong
              \tau^{-1}H(1)_\lambda$ for some $\lambda\in\mathbb{P}_k^1$.
  \end{itemize}
\end{lemm}

\begin{proof}(1) and (2) are obvious.

(3) Note that $M=X[2]$ for a quasi-simple module $X$ with dimension
vector $(1,1,1)$.  Thus $\{1,2,3\}<\mu(M)<\mu(I_3)=\{1,2,3,4,5\}$,
the maximal GR measure.  Thus $\mu(M)=\{1,2,3,5\}$.

(4) $M=X[2]$ for a quasi-simple module $X$ with dimension vector
$(0,1,1)$.  It is easily seen that $\Hom(H^1,M)=0=\Hom(\tau
H(1)_\lambda,M)$ for each $\lambda\in\mathbb{P}_k^1$. Thus
$\{1,2,8\}\leq \mu(M)<\{1,2,3\}$. In particular the GR submodule $Y$
of $M$ is regular and thus $\tau Y$ is a submodule of $\tau M$ with
$\udim \tau M=(1,2,2)$. It is not difficult to see that the only
possibility is $\udim Y=(0,1,1)$. Thus $\mu(M)=\{1,2,8\}$.

(5) If $\udim M=(1,2,3)$, then $M\cong
              \tau^{-1}H(1)_\lambda$ for some $\lambda\in\mathbb{P}_k^1$.
The GR submodule of $M$ is the projective module $P_1$. Conversely,
if $\mu(M)=\{1,3,4,6\}$, then $M$ contains $P_1$ as a GR submodule
and the corresponding GR factor has dimension vector $(0,1,1)$. Thus
$\udim M=\{1,2,3\}$.
\end{proof}

We have seen that an indecomposable regular $M$ with $(\udim M)_1=1$
has quasi-length at most 2 and $\ql(M)=2$ implies that $\udim
M=(1,2,2)$ or $(1,3,4)$, i.e., $M$ is in the AR component containing
$H(1)_\lambda$ for some $\lambda\in\mathbb{P}_k^1$.

\begin{lemm}
   Let $M$ be an indecomposable module with dimension
             $(1,2,2)$ or $(1,3,4)$  and $\ql(M)=2$. If  $N$ is a quasi-simple
             module  with $\udim N=\udim M$, then $\mu(N)<\mu(M)$.
\end{lemm}

\begin{proof}Let $\udim M=(1,2,2)$. We have seen that  $\mu(M)=\{1,2,3,5\}$.
If $N$ is a quasi-simple module with $\udim N=(1,2,2)$, then
$\mu(N)$ does not start with $\{1,2,3\}$. Otherwise, $N$ contains
some indecomposable module $X$ with dimension vector $(1,1,1)$ or
$(0,2,1)$ as a GR submodule. If $\udim X=(1,1,1)$, then there is an
epimorphism from $X[2]$ to $N$ which is impossible since $\udim
X[2]=(1,2,2)=\udim N$. Note that $\udim X=(0,2,1)$ is either not
possible since otherwise the factor contains the projective simple
module.  Therefore, $\mu(N)<\{1,2,3\}<\mu(M)$.

Now let $\udim M=(1,3,4)$ and $N$ be a quasi-simple module with
dimension vector $\udim N=\udim M$. Then $N$ does not contain any
$H(1)_\lambda$ as a submodule since otherwise $\tau N$, with $\udim
\tau N=(1,2,2)$, contains a submodule module with dimension vector
$(1,1,1)$. This is not possible by above discussion. Then we again
have $\mu(N)<\{1,2\}<\mu(M)$.
\end{proof}

We may ask the following question in general:
\begin{Question} Let
$M$ and $N$ be indecomposable regular modules with dimension vector
$\udim M=\udim N$ and $\ql(M)>\ql(N)$. Does $\mu(M)>\mu(N)$ hold?
\end{Question}

Namely, we may calculate precisely the GR measures of the
quasi-simple modules $N$ with $\udim N=(1,3,4)$ and $M$ with $\udim
M=(1,2,2)$. Since $\mu(N)<\{1,2\}$, $\mu(N)$ starts with
$\mu(P_1)=\{1,3,4\}$. On the other hand, $\mu(N)$ does not contain
$7$. Namely, assume that it is not the case and let $X$ be a GR
submodule of $N$. Then the GR factor module $N/X$ is simple and thus
$\udim X=(0,3,4)$ or $(1,2,4)$. However, the first vector correspond
the preprojective module $H_3$ and the second vector is not a root.
Similarly, a detailed discussion shows that $\mu(N)$ does not
contains $5$. Therefore, the only possibility is that
$\mu(N)=\{1,3,4,6,8\}$. Thus any GR submodule of $N$ is of dimension
vector $(1,2,3)$.  It follows that any quasi-simple module $M$ with
$\udim M=(1,2,2)$ contains a submodule of dimension vector
$(1,2,3)\Phi=(0,1,1)$. Thus the  GR measure of $M$ is either
$\{1,2,4,5\}$ or $\{1,2,5\}$, and both possibilities occur.

\smallskip
We have known the GR measures of the indecomposable modules $M$ with
lengths not larger than $6$ except $\udim M=(2,2,2)$ and $(1,3,2)$.
For $\udim M=(2,2,2)$, we have seen that the GR factor module is the
simple injective. Thus the GR submodule of $M$ has dimension vector
$(1,2,2)$. Note that $M$ always contains an indecomposable module
with dimension vector $(1,1,1)$.  Therefore, $\mu(M)=\{1,2,3,5,6\}$.

There are several possibilities for  $\udim(M)=(1,3,2)$. Note that
the length of any GR factor $Y$ of $M$ does not equal to $2$. Since
otherwise, the GR submodule $X$ of $M$ has dimension vector
$(1,2,1)$. But there does not exist an epimorphism from $\tau^{-1}X$
to $Y$ because $\udim \tau^{-1}X=(1,1,1)$ and $\udim Y=(0,1,1)$
(Proposition \ref{moepi}(4)). Thus $Y$ is a simple module or the
length of $Y$ is $3$. Thus only the following possibilities may
occur: $\{1,2,3,6\}$, $\{1,2,3,5,6\}$, $\{1,2,4,5,6\}$,
$\{1,2,5,6\}$.

\section{Gabriel-Roiter measures admitting no direct predecessors}

Ringel showed in \cite{R2} that each GR measure different from
$\mu^1$, the maximal GR measure, has a direct successor. However,
there are GR measures admitting no direct predecessors in general. A
GR measure for an algebra of finite representation type definitely
has a direct predecessor.  Moreover,  we have shown in \cite{Ch6}
that for path algebras $kQ$ of tame quivers, $1\leq
\ndp(\Lambda)<\infty$ (meaning that there are only finitely many GR
measures processing no direct predecessors). Thus it is natural to
ask if the number of the GR measures admitting no direct
predecessors relates to the representation type (of hereditary
algebras). More precisely, we want to know whether an algebra
$\Lambda$ (over an algebraically closed field) is of wild type
implies that there are infinitely many GR measures admitting no
direct predecessors, i.e., $\ndp(\Lambda)=\infty$, and vice versa.

It should be very difficult to answer this question in general. Now
we come back to the quiver $Q$ we considered above. We have seen
$\mu(P_1)=\{1,3,4\}$ is the minimal central measure and thus has no
direct predecessor. In this section,  we will show the following
theorem:

\begin{theo}\label{ndp} Let $\mu^n=\{1,2,4,\ldots,2n,2n+1\}$, $n\geq 1$.
Then $\mu^n$ has no direct predecessor.
\end{theo}

\begin{lemm}
  \begin{itemize}
     \item[(1)] $\mu^n$ is a GR measure.
     \item[(2)] If $M$ is an indecomposable module with $\mu(M)=\mu^n$. Then $\udim M=
                (1,n,n)$ or $(0,n+1,n)$.
   \end{itemize}
\end{lemm}

\begin{proof} It is known that each indecomposable module with dimension vector $(0,n+1,n)$
has GR measure $\mu^n$. Thus $\mu^n$ is a GR measure. We have seen
that a non-injective indecomposable module $M$ has GR measure
$\{1,2,4,\ldots,2n\}$ if and only if $\udim M=(0,n,n)$.  Thus an
indecomposable module with GR measure $\mu^n$ has dimension vector
$(1,n,n)$ or $(0,n+1,n)$.
\end{proof}

\begin{lemm}\label{a} Let $M$ be an indecomposable module with GR measure $\mu(M)=\mu^n$. Then
        each indecomposable regular factor module of $M$ contains some indecomposable module with
    dimension vector $(0,1,1)$ as a submodule.
\end{lemm}

\begin{proof} By above lemma, $\udim M= (1,n,n)$ or $(0,n+1,n)$. In each case, we have a
short exact sequence
$$0\ra X\stackrel{\iota}{\ra} M\ra M/X\ra 0$$
where  $\iota$ is a GR inclusion and thus $\udim X=(0,n,n)$. Note
that the factor $M/X$ is a preinjective simple module.  Let
$M\stackrel{\pi}{\ra}Y$ be an epimorphism with $Y$ an indecomposable
regular module.  Then $\Hom(M/X,Y)=0$. By definition of cokernel,
the composition $X\stackrel{\pi\iota}{\ra}Y$ is not zero. Since an
indecomposable non-simple factor of $X$ has dimension vector
$(0,a,a)$ or $(0,a+1,a)$,  the image of $\pi\iota$  contains a
submodule with dimension vector $(0,1,1)$.
\end{proof}

\begin{lemm}\label{b} Fix an $n\geq 1$. Let $M$ be an indecomposable module such that
$\mu^n<\mu(M)$.  Then $\mu(M)$ starts with
$\mu(M)=\{1,2,\ldots,2m,2m+1\}$ for some $1\leq m\leq n$. In
particular, $M$ contains an indecomposable submodule with GR measure
$\mu^m$.
\end{lemm}

\begin{proof}
This follows from the definition of GR measure.
\end{proof}

\begin{lemm}\label{c} If $M$ is an indecomposable module such that
$\mu=\mu(M)$ is a direct predecessor of $\mu^n$ for some $n$. Then
$M$ is regular.
\end{lemm}

\begin{proof}
By the calculation of the GR submodules, we have $\mu(P)<\{1,2\}$
for every indecomposable preprojective module $P$. Note that for
each indecomposable preprojective module, there are infinitely many
preprojective modules with greater GR measures. Thus $M$ is not
preprojective. On the other hand, let $X$ be an indecomposable
regular with dimension vector $(1,1,1)$ and $I$ be an indecomposable
preinjective module such that $I\ncong S_2$. Then $\Hom(X,I)\neq 0$.
Therefore, if $I$ is neither isomorphic to the simple modules $S_1$,
$S_2$,  nor isomorphic to the injective module $I_2$, there is
always a monomorphism from $X$ to $I$ and thus $\mu(I)$ starts with
$\{1,2,3\}>\mu^n$ for every $n\geq 1$. It follows that $M$ has to be
a regular module.
\end{proof}

{\it proof of Theorem.} For the purpose of a contradiction, we
assume that $M$ is an indecomposable module such that $\mu(M)$ is a
direct predecessor of $\mu^n$ for a fixed $n$. Thus by Lemma
\ref{c}, we may write $M=X[r]$ for some quasi-simple $X$ and $r\geq
1$. It follows that $\mu(M)=\mu(X[r])<\mu^n<\mu(X[r+1])$. Thus
$X[r+1]$ contains a submodule $Y$ with GR measure $\mu_m$ for some
$1\leq m\leq n$ (Lemma \ref{b}). Note that $\udim Y=(1,m,m)$ or
$(0,m+1,m)$ and $\mu(Y)\geq \mu^n$. We claim that
$\Hom(Y,\tau^{-r}X)=0$. If this is not the case, then by Lemma
\ref{a}, the image of a nonzero homomorphism, in particular
$\tau^{-r}X$, contains a submodule $Z$ with dimension vector
$(0,1,1)$. Therefore, there is a monomorphism $\tau^r Z\ra X$, and
thus
 $$\mu^n>\mu(M)\geq \mu(X)>\mu(\tau^r Z)\geq \{1,2,3\},$$ which is a contradiction.
Since there is a short exact sequence $$0\ra M=X[r]\ra X[r+1]\ra
\tau^{-r}X\ra 0$$ and $\Hom(Y,\tau^{-r}X)=0$, the inclusion $Y\ra
X[r+1]$ factors through $X[r]$. Therefore, there is a monomorphism
$Y\ra X[r]$. It follows that
$$\mu(X[r])\geq \mu(Y)\geq \mu^n=\mu(M).$$  This contradiction implies that $\mu^n$ has no direct
predecessor for each $n\geq 1$. \qed

\end{document}